\documentclass[a4paper, 12pt]{article}
\usepackage{amssymb, amsmath, amsfonts}
\usepackage{amsthm,amscd,amsfonts,latexsym,color,epsfig}
\begin{document}
	
	\begin{center}
		{\bf Unique solvability of a non-local problem for mixed type equation with fractional derivative}\\
		 E.T. Karimov, A.S.Berdyshev, Rakhmatullaeva N.A.\\
		E-mail: erkinjon@gmail.com, berdyshev@mail.ru, rakhmatullaeva@mail.ru\\
		\emph{Institute of Mathematics, National University of Uzbekistan (Tashkent, Uzbekistan),}\\
		\emph{Kazakh National Pedagogical University (Almaty, Kazakhstan),}\\
			\emph{Tashkent State Technical University (Tashkent, Uzbekistan)}\\
	\end{center}
	
	\bigskip
	
	\textbf{MSC 2000:} 35M10\\
	\textbf{Keywords:} Caputo fractional derivative; mixed type equation; Volterra integral equation; parabolic-hyperbolic type equation; Green's function
	
	\bigskip
	
	\textbf{Abstract. }In this work we investigate a boundary problem with non-local conditions for mixed parabolic-hyperbolic type equation with three lines of type-changing with Caputo fractional derivative in the parabolic part. We equivalently reduce considered problem to the system of second kind Volterra integral equations. In the parabolic part we use solution of the first boundary problem with appropriate Green's function and in hyperbolic parts we use corresponding solutions of the Cauchy problem.  
	
	\bigskip

 \section{Formulation of a problem}

Due to many applications in gas and aerodynamics, mechanics, theory of boundary problems for mixed type equations was very rapidly developed [1]. Nowadays, this theory has many branches due to the usage of various methods of mathematical and functional analysis, topological methods and the method of Fractional Calculus [2-4].

We note closely-related recent works on local and non-local boundary problems for parabolic-hyperbolic equations [5-9].

Regarding the investigations of mixed parabolic-hyperbolic equations with non-smooth lines of type changing we note works [10-14].

In this work we investigate boundary problem with non-local condition, which connects value of seeking function in one characteristics with value of this function on another characteristics. This kind of non-local condition was used in the work [15]. We study this problem for unique solvability. Spectral properties of such problems in integer case were investigated in works [16-17].

Consider an equation
$$
0=\left\{
\begin{array}{l}
u_{xx}-_CD_{0y}^\lambda u,\,\,\,(x,y)\in\Omega_0,\\
u_{xx}-u_{yy},\,\,\,\,(x,y)\in\Omega_i(i=\overline{1,3})\\
\end{array}
\right.
\eqno (1.1)
$$
in a domain $\Omega=\Omega_0\cup \Omega_1\cup \Omega_2\cup \Omega_3\cup AB\cup AA_0\cup AB_0$. Here $\Omega_0=\left\{(x,y):\, 0<x<1,\,0<y<1\right\}$, $\Omega_1,\,\Omega_2,\,\Omega_3$ are characteristic triangles with endpoints $A,B,C;\, A,A_0,D;\, B,B_0,E$ respectively, where $A(0,0), A_0(0,1)$, $B(1,0), B_0(1,1)$, $C(1/2,-1/2), D(-1/2,1/2)$, $E(3/2,1/2)$, $0<\lambda\leq 1$,
$$
_CD_{0y}^\lambda g =\left\{\begin{array}{l}
\frac{dg}{dy},\hfill \lambda=1\\
\frac{1}{\Gamma(1-\lambda)}\int\limits_0^y (y-t)^{-\lambda} g'(t)dt,\hfill 0<\lambda<1
\end{array}\right.
$$
is the Caputo fractional operator [4].

\textbf{Problem. } Find a solution of the equation (1.1) from the class of functions 
$$
W=\left\{u(x,y):\, u\in C(\overline{\Omega})\cap C^1(\Omega\ AB)\cap C_{x,y}^{2,\lambda}(\Omega_0)\cap C^2(\Omega_i)\right\},
$$
satisfying conditions
$$
a_1(t)u(t,-t)+a_2(t)u(t,-t)=a_3(t),\,\,0\leq t\leq 1, \eqno (1.2)
$$
$$
u|_{CB}=\varphi_1(x),\,\,\,1/2\leq x\leq 1, \eqno (1.3)
$$
$$
u|_{BE}=\varphi_2(y),\,\,\, 0\leq y\leq 1/2, \eqno (1.4)
$$
$$
\lim\limits_{y\rightarrow +0}{\mathop{\lim }}\,{}_{C}D_{0y}^\lambda u(x,y)=u_y(x,-0),\,\,\,0<x<1, \eqno (1.5)
$$
where $a_i(\cdot),\, \varphi_j(\cdot) (i-\overline{1,3}, j=1,2)$ are given functions such that $\varphi_1(1)=\varphi_2(0)$.

\section{Main functional correlations}

Solution of the Cauchy problem for the equation (1.1) in the domain $\Omega_1$ can be represented as [18]
$$
u(x,y)=\frac{1}{2}\left[\tau_1(x+y)+\tau_1(x-y)+\int\limits_{x-y}^{x+y}\nu_1^-(t)dt,\right]\eqno (2.1)
$$
where $u(x,0)=\tau_1(x),\,u_y(x,-0)=\nu_1^-(x)$.

Solution (2.1) satisfy to the condition (1.3) and after some evaluations we get
$$
\nu_1^{-}(x)=-\tau_1'(x)+\varphi_1\left(\frac{x+1}{2}\right),\,\,\,0\leq x\leq 1. \eqno (2.2)
$$
Equality (2.2) is one the main functional correlations on $AB$ reduced from the hyperbolic part of the mixed domain. Now let us deduce another functional correlations between functions $\tau_1(x)$,  $\nu_1(x)$. For this aim we pass to the limit as $y\rightarrow +0$ and from the equation (1.1) we obtain [4]
$$
\tau_1''(x)-\Gamma(\lambda)\nu_1^{+}(x)=0, \eqno (2.3)
$$
where $\nu_1^{+}(x)=\lim\limits_{y\rightarrow +0} {\mathop{\lim }}\,{}_{C}D_{0y}^\lambda u(x,y)$.

From the equalities (2.2) and (2.3), considering transmitting condition (1.5) we get
$$
\tau_1''(x)+\Gamma(\lambda)\tau_1'(x)=\varphi_1'\left(\frac{x+1}{2}\right). \eqno (2.4)
$$
General solution of (2.4) has a form
$$
\tau_1(x)=C_1+e^{-\Gamma(\lambda)x}\left[C_2+2\int\limits_0^x e^{\Gamma(\lambda)t}\varphi_1'\left(\frac{t+1}{2}\right)dt\right].
$$

From the condition (1.2) and (1.3) we have
$$
\tau_1(0)=\frac{a_3(0)}{a_1(0)+a_2(0)},\,\,\tau_1(1)=\varphi_1(1). \eqno (2.5)
$$

Solution of the equation (2.3), satisfying conditions (2.5) we represent as 
$$
\begin{array}{l}
\displaystyle{\tau_1(x)=2\int\limits_0^x e^{\Gamma(\lambda)(t-x)}\varphi_1\left(\frac{t+1}{2}\right)dt+\frac{e^{\Gamma(\lambda)(1-x)}-1}{e^{\Gamma(\lambda)}-1}\frac{a_3(0)}{a_1(0)+a_2(0)}+
\frac{e^{\Gamma(\lambda)}-e^{\Gamma(\lambda)(1-x)}}{e^{\Gamma(\lambda)}-1}\varphi_1(1)+}\\
+\displaystyle{\frac{2}{e^{\Gamma(\lambda)}-1}\int\limits_0^1 
\left[e^{\Gamma(\lambda)(t-x)}-e^{\Gamma(\lambda)t}\right]\varphi_1\left(\frac{t+1}{2}\right)dt.}\\
\end{array}
\eqno (2.6)
$$
Formula (2.6) gives explicit representation of the function $\tau_1(x)$. Using (2.2) and (1.5) we find functions $\nu_1^-(x), \, \nu_1^+(x)$ in an explicit form.

Using (2.1) we find $u(t,-t)$:
$$
u(t,-t)=\frac{1}{2}\left[\tau_1(0)+\tau_1(2t)-\int\limits_0^{2t}\nu_1^-(z)dz\right],\,\,0\leq t\leq 1/2. \eqno (2.7)
$$

Solution of the Cauchy problem for (1.1) in the domain $\Omega_2$ has a form
$$
u(x,y)=\frac{1}{2}\left[\tau_2(y+x)+\tau_2(y-x)+\int\limits_{y-x}^{y+x}\nu_2(t)dt,\right],\eqno (2.8)
$$
where $u(0,y)=\tau_2(y) (0\leq y\leq 1),\,u_x(0,y)=\nu_2(y) (0<y<1)$.

Using this representation calculate $u(-t,t)$:
$$
u(-t,t)=\frac{1}{2}\left[\tau_2(0)+\tau_2(2t)-\int\limits_0^{2t}\nu_2(z)dz\right],\,\,0\leq t\leq 1/2. \eqno (2.9)
$$

Equalities (2.7) and (2.9) we substitute into (1.2) and after differentiation once with respect to $t$, we obtain
$$
\begin{array}{l}
\displaystyle{\nu_2(t)=\tau_2'(t)+\frac{a_2(t/2)}{a_1(t/2)}\left[\tau_1'(t)-\nu_1(t)\right]+\frac{a_1(t/2)a_2'(t/2)-a_1'(t/2)a_2(t/2)}{2a_1^2(t/2)}\times}\\
\times \displaystyle{\left[\tau_1(0)+\tau_1(t)-\int\limits_0^t\nu_1(z)dz\right]+\frac{a_1'(t/2)a_3(t/2)-a_3'(t/2)a_1(t/2)}{a_1^2(t/2)}.}\\
\end{array}
\eqno (2.10)
$$

Now we write solution of the Cauchy problem for (1.1) in the domain $\Omega_2$:
$$
u(x,y)=\frac{1}{2}\left[\tau_3(y+x-1)+\tau_3(y-x+1)+\int\limits_{y+x-1}^{y-x+1}\nu_3(t)dt,\right],\eqno (2.11)
$$
where $u(1,y)=\tau_3(y) (0\leq y\leq 1),\,u_x(1,y)=\nu_3(y) (0<y<1)$.

(2.11) satisfy to (1.4) and differentiating once we get
$$
\nu_3(y)=-\tau_3'(y)+\varphi_2\left(\frac{y}{2}\right),\,\,0<y<1. \eqno (2.12)
$$

Equalities (2.10) and (2.12) are main functional correlations, which we shall use further. 

\section{Reduction to the system of integral equations}

Solution of the first boundary problem for the equation (1.1) in the domain ${{\Omega }_{0}}$ has a form [4]
$$
  u\left( x,y \right)=\int\limits_{0}^{y}{{{G}_{{{x}_{1}}}}\left( x,y;0,{{y}_{1}} \right){{\tau }_{2}}\left( {{y}_{1}} \right)d{{y}_{1}}}-\int\limits_{0}^{y}{{{G}_{{{x}_{1}}}}\left( x,y;1,{{y}_{1}} \right){{\tau }_{3}}\left( {{y}_{1}} \right)d{{y}_{1}}}+
\int\limits_{0}^{1}{\overline{G}\left( x-{{x}_{1}},y\right)\tau_1(x_1)dx_1,}
\eqno (3.1)
$$
where
$$
\overline{G}\left( x-{{x}_{1}},y\right)=\frac{1}{\Gamma(1-\lambda)}\int\limits_0^yy_1^{-\lambda}G\left( x,y;{{x}_{1}},{{y}_{1}} \right)dy_1,
$$
$$
G\left( x,y;{{x}_{1}},{{y}_{1}} \right)=\frac{{{\left( y-{{y}_{1}} \right)}^{\rho -1}}}{2}\sum\limits_{n=0}^{\infty }{\left[ e_{1,\rho }^{1,\rho }\left( -\frac{\left| x-{{x}_{1}}+2n \right|}{{{\left( y-{{y}_{1}} \right)}^{\rho }}} \right)-e_{1,\rho }^{1,\rho }\left( -\frac{\left| x+{{x}_{1}}+2n \right|}{{{\left( y-{{y}_{1}} \right)}^{\rho }}} \right) \right]}
$$
is the Green's function of the first boundary problem [4], $\rho =\lambda /2$,
$$
e_{\alpha ,\beta }^{\mu ,\delta }\left( z \right)=\sum\limits_{n=1}^{\infty }{\frac{{{z}^{n}}}{\Gamma \left( \alpha n+\mu  \right)\Gamma \left( \delta -\beta n \right)}}
$$
is the Wright type function [4].

Using the following formulas [4]
$$
\begin{array}{l}
\displaystyle{
\frac{d}{dz}e_{\alpha ,\rho }^{\mu ,\delta }\left( z \right)=-\frac{1}{\rho z}\left[ e_{\alpha ,\rho }^{\mu ,\delta -1}\left( z \right)+\left( 1-\delta  \right)e_{\alpha ,\rho }^{\mu ,\delta }\left( z \right) \right],}\\
\displaystyle{e_{1,\rho }^{1,\delta -1}\left( z \right)+\left( 1-\delta  \right)e_{1,\rho }^{1,\delta }\left( z \right)=\rho ze_{1,\rho }^{1,\delta -\rho }\left( z \right),}\\
\displaystyle{
\frac{d}{dz}\left( {{z}^{\mu -1}}e_{\alpha ,\rho }^{\mu ,\delta }\left( c{{z}^{\alpha }} \right) \right)={{z}^{\mu -2}}e_{\alpha ,\rho }^{\mu -1,\delta }\left( c{{z}^{\alpha }} \right),}\\
\displaystyle{
 \frac{d}{dz}\left( {{z}^{\mu -1}}e_{\alpha ,\rho }^{\mu ,\delta }\left( c{{z}^{-\rho }} \right) \right)={{z}^{\delta -\alpha }}e_{\alpha ,\rho }^{\mu ,\delta -1}\left( c{{z}^{-\rho }} \right),}\\
\end{array}
\eqno (3.2)
$$
after some evaluations we deduce
$$
{{G}_{{{x}_{1}}x}}\left( +0,y;0,{{y}_{1}} \right)=\frac{d}{d{{y}_{1}}}\left( \sum\limits_{n=-\infty }^{\infty }{\frac{1}{{{\left( y-{{y}_{1}} \right)}^{\rho }}}e_{1,\rho }^{1,1-\rho }\left( -\frac{\left| 2n \right|}{{{\left( y-{{y}_{1}} \right)}^{\rho }}} \right)} \right),\eqno (3.3)
$$
$$
{{G}_{{{x}_{1}}x}}\left( +0,y;1,{{y}_{1}} \right)=\frac{d}{d{{y}_{1}}}\left( \sum\limits_{n=-\infty }^{\infty }{\frac{1}{{{\left( y-{{y}_{1}} \right)}^{\rho }}}e_{1,\rho }^{1,1-\rho }\left( -\frac{\left| 2n+1 \right|}{{{\left( y-{{y}_{1}} \right)}^{\rho }}} \right)} \right),\eqno (3.4)
$$
$$
{{G}_{{{x}_{1}}x}}\left( 1-0,y;0,{{y}_{1}} \right)=\frac{d}{d{{y}_{1}}}\left( \sum\limits_{n=-\infty }^{\infty }{\frac{1}{{{\left( y-{{y}_{1}} \right)}^{\rho }}}e_{1,\rho }^{1,1-\rho }\left( -\frac{\left| 2n+1 \right|}{{{\left( y-{{y}_{1}} \right)}^{\rho }}} \right)} \right),\eqno (3.5)
$$
$$
{{G}_{{{x}_{1}}x}}\left( 1-0,y;1,{{y}_{1}} \right)=\frac{d}{d{{y}_{1}}}\left( \sum\limits_{n=-\infty }^{\infty }{\frac{1}{{{\left( y-{{y}_{1}} \right)}^{\rho }}}e_{1,\rho }^{1,1-\rho }\left( -\frac{\left| 2n \right|}{{{\left( y-{{y}_{1}} \right)}^{\rho }}} \right)} \right).\eqno (3.6)
$$

Differentiating (3.1) once with respect to $x$ and at $x\rightarrow +0$, $x\rightarrow 1-0$, considering (3.3)-(3.6), we obtain
$$
\begin{array}{l}
\displaystyle{
  u_x\left( +0,y \right)=\nu _{2}\left( y \right)=\int\limits_{0}^{y}{{{\tau }_{2}}\left( {{y}_{1}} \right)\frac{\partial }{\partial {{y}_{1}}}\left( \sum\limits_{n=-\infty }^{\infty }{\frac{e_{1,\rho }^{1,1-\rho }\left( -\frac{\left| 2n \right|}{{{\left( y-{{y}_{1}} \right)}^{\rho }}} \right)}{{{\left( y-{{y}_{1}} \right)}^{\rho }}}} \right)d{{y}_{1}}}-} \\
\displaystyle{  -\int\limits_{0}^{y}{{{\tau }_{3}}\left( {{y}_{1}} \right)\frac{\partial }{\partial {{y}_{1}}}\left( \sum\limits_{n=-\infty }^{\infty }{\frac{e_{1,\rho }^{1,1-\rho }\left( -\frac{\left| 2n+1 \right|}{{{\left( y-{{y}_{1}} \right)}^{\rho }}} \right)}{{{\left( y-{{y}_{1}} \right)}^{\rho }}}} \right)d{{y}_{1}}}+\int\limits_{0}^{1}{\overline{G}_x\left(-{{x}_{1}},y\right)\tau_1(x_1)dx_1}} \\
\end{array}\eqno (3.7)
$$
$$
\begin{array}{l}
\displaystyle{
  {{u}_{x}}\left( 1-0,y \right)=\nu _{3}\left( y \right)=\int\limits_{0}^{y}{{{\tau }_{2}}\left( {{y}_{1}} \right)\frac{\partial }{\partial {{y}_{1}}}\left( \sum\limits_{n=-\infty }^{\infty }{\frac{e_{1,\rho }^{1,1-\rho }\left( -\frac{\left| 2n+1 \right|}{{{\left( y-{{y}_{1}} \right)}^{\rho }}} \right)}{{{\left( y-{{y}_{1}} \right)}^{\rho }}}} \right)d{{y}_{1}}}-} \\
\displaystyle{  -\int\limits_{0}^{y}{{{\tau }_{3}}\left( {{y}_{1}} \right)\frac{\partial }{\partial {{y}_{1}}}\left( \sum\limits_{n=-\infty }^{\infty }{\frac{e_{1,\rho }^{1,1-\rho }\left( -\frac{\left| 2n \right|}{{{\left( y-{{y}_{1}} \right)}^{\rho }}} \right)}{{{\left( y-{{y}_{1}} \right)}^{\rho }}}} \right)d{{y}_{1}}}+\int\limits_{0}^{1}{\overline{G}_x\left(1-{{x}_{1}},y\right)\tau_1(x_1)dx_1}.} \\
\end{array}\eqno (3.8)
$$

Using formula of integration by parts, (3.7) and (3.8) we rewrite as follows
$$
\nu _{2}\left( y \right)=\int\limits_{0}^{y}{{{\tau }_{2}}'\left( {{y}_{1}} \right){{K}_{1}}\left( y,{{y}_{1}} \right)d{{y}_{1}}}-\int\limits_{0}^{y}{{{\tau }_{3}}'\left( {{y}_{1}} \right){{K}_{2}}\left( y,{{y}_{1}} \right)d{{y}_{1}}}+\int\limits_{0}^{1}{\overline{G}_x\left(-{{x}_{1}},y\right)\tau_1(x_1)dx_1},
\eqno (3.9)
$$
$$
\nu _{3}\left( y \right)=\int\limits_{0}^{y}{{{\tau }_{2}}'\left( {{y}_{1}} \right){{K}_{2}}\left( y,{{y}_{1}} \right)d{{y}_{1}}}-\int\limits_{0}^{y}{{{\tau }_{3}}'\left( {{y}_{1}} \right){{K}_{1}}\left( y,{{y}_{1}} \right)d{{y}_{1}}}
+\int\limits_{0}^{1}{\overline{G}_x\left(1-{{x}_{1}},y\right)\tau_1(x_1)dx_1},
\eqno (3.10)
$$
where
$$
{{K}_{1}}\left( y,{{y}_{1}} \right)=\frac{{{\left( y-{{y}_{1}} \right)}^{-\rho }}}{\Gamma \left( 1-\rho  \right)}+\sum\limits_{\begin{smallmatrix}
 n=-\infty  \\
 n\ne 0
\end{smallmatrix}}^{\infty }{\frac{1}{{{\left( y-{{y}_{1}} \right)}^{\rho }}}}e_{1,\rho }^{1,1-\rho }\left( -\frac{\left| 2n \right|}{{{\left( y-{{y}_{1}} \right)}^{\rho }}} \right)
=\frac{{{\left( y-{{y}_{1}} \right)}^{-\rho }}}{\Gamma \left( 1-\rho  \right)}+{{\widetilde{K}}_{1}}\left( y,{{y}_{1}} \right),
\eqno     (3.11)
$$
$$
{{K}_{2}}\left( y,{{y}_{1}} \right)=\sum\limits_{n=-\infty }^{\infty }{\frac{1}{{{\left( y-{{y}_{1}} \right)}^{\rho }}}}e_{1,\rho }^{1,1-\rho }\left( -\frac{\left| 2n+1 \right|}{{{\left( y-{{y}_{1}} \right)}^{\rho }}} \right).\eqno (3.12)
$$

Considering (2.10), (2.12) and (3.9), (3.10) we get the following system of second kind Volterra integral equations:
$$
\left\{
\begin{array}{l}
\displaystyle{\tau_2(y)-\int\limits_{0}^{y}{{{\tau }_{2}}'\left( {{y}_{1}} \right){{K}_{1}}\left( y,{{y}_{1}} \right)d{{y}_{1}}}=-
\int\limits_{0}^{y}{{{\tau }_{3}}'\left( {{y}_{1}} \right){{K}_{2}}\left( y,{{y}_{1}} \right)d{{y}_{1}}}+f_1(y),}\\
\displaystyle{\tau_3(y)-\int\limits_{0}^{y}{{{\tau }_{3}}'\left( {{y}_{1}} \right){{K}_{1}}\left( y,{{y}_{1}} \right)d{{y}_{1}}}=-
\int\limits_{0}^{y}{{{\tau }_{2}}'\left( {{y}_{1}} \right){{K}_{2}}\left( y,{{y}_{1}} \right)d{{y}_{1}}}+f_2(y),}\\
\end{array}
\right.\eqno (3.13)
$$
where
$$
\begin{array}{l}
\displaystyle{f_1(y)=\int\limits_{0}^{1}{\overline{G}_x\left(-{{x}_{1}},y\right)\tau_1(x_1)dx_1}+
\frac{a_2\left(\frac{t}{2}\right)}{a_1\left(\frac{t}{2}\right)}\left[\nu_1(t)-\tau_1'(t)\right]+\frac{a_1\left(\frac{t}{2}\right)a_2'\left(\frac{t}{2}\right)
-a_1'\left(\frac{t}{2}\right)a_2\left(\frac{t}{2}\right)}{2a_1^2\left(\frac{t}{2}\right)}\times}\\
\times \displaystyle{\left[\int\limits_0^t\nu_1(z)dz-\tau_1(0)-\tau_1(t)\right]-\frac{a_1'\left(\frac{t}{2}\right)a_3\left(\frac{t}{2}\right)-a_3'\left(\frac{t}{2}
\right)a_1\left(\frac{t}{2}\right)}{a_1^2\left(\frac{t}{2}\right)},}\\
\displaystyle{f_2(y)=\varphi_2\left(\frac{y}{2}\right)-\int\limits_{0}^{1}{\overline{G}_x\left(1-{{x}_{1}},y\right)\tau_1(x_1)dx_1}}\\
\end{array}
\eqno (3.14)
$$
Considering representations of kernels (3.11), (3.12) and right-hand sides (3.14), one can easily be convinced that system of integral equations (3.13) is uniquely solvable. After we find functions $\tau_j(y),\,\nu_j(y) (j=2,3)$, we recover solution of the formulated problem in the domain $\Omega_0$ by the formula (3.1) and in the domains $\Omega_1,\,\Omega_2,\,\Omega_3$ by formulas (2.1), (2.8), (2.11), respectively.

Finally we can formulate our result as the following theorem:

\textbf{Theorem.} If $a_i(\cdot),\, \varphi_2(\cdot)\in C^1[0,1/2]\cap C^2(0,1/2)$, $\varphi_1(\cdot)\in C^1[1/2,1]\cap C^2(1/2,1)$,
then there exist unique solution of the problem (1.1)-(1.5)

\section{Acknowledgement}
This research was partially supported by the Grant No 3293GF4 of the Ministry of education and science of the Republic of Kazakhstan

\end{document}